\documentclass[letterpaper, 10 pt, conference, english]{ieeeconf} 
\IEEEoverridecommandlockouts
\usepackage{epsfig}
\usepackage{subfigure}
\usepackage{amssymb,amsmath,amsfonts,layout,graphicx}
\usepackage{makeidx}
\usepackage{babel}
\usepackage{sublabel}
\usepackage{cases}
\usepackage{algorithm}
\usepackage{algorithmic}

 \usepackage{pgfplots}
  \pgfplotsset{compat=newest}
  \usetikzlibrary{plotmarks}
  \usetikzlibrary{arrows.meta}
  \usepgfplotslibrary{patchplots}
  \usepackage{grffile}
  \usepackage{amsmath}

\usepackage
[
        letterpaper,
        left=1.91cm,
        right=1.91cm,
        top=1.91cm,
        bottom=2cm,
]
{geometry}

\newtheorem{proposition}{Proposition}

 \title{Off-Street Parking for TNC Vehicles to Reduce Cruising Traffic}

\author{Sen Li$^\star$, Junjie Qin$^\dagger$, Hai Yang$^\star$, Kameshwar Poolla$^\dagger$ and Pravin Varaiya$^\ddagger$
\thanks{$^\star$The Department of Civil and Environmental Engineering, The Hong Kong University of Science and Technology, Clear Water Bay, Hong Kong. Email: \{cesli, cehyang\}@ust.hk}
\thanks{$^\dagger$The Department of Mechanical Engineering, University of California, Berkeley, USA. Email: \{qinj, poolla\}@berkeley.edu}
 \thanks{$^\ddagger$The Department of Electrical Engineering and Computer Science, University of California, Berkeley, USA. Email: varaiya@berkeley.edu}
}

\begin{document}
\maketitle
\begin{abstract} 
This paper considers off-street parking for the cruising vehicles of transportation network companies (TNCs) to reduce the traffic congestion. We propose a novel business that integrates the shared parking service into the TNC platform. In the proposed model, the platform (a) provides interfaces that connect passengers, drivers and garage operators (commercial or private garages); (b) determines the ride fare, driver payment, and parking rates; (c) matches passengers  to TNC vehicles for ride-hailing services; and (d) matches vacant TNC vehicles to unoccupied parking garages to reduce the cruising cost. A queuing-theoretic model is proposed to capture the matching process of passengers, drivers, and parking garages.  A market-equilibrium model is developed to capture the incentives of the passengers, drivers, and garage operators. An  optimization-based model is formulated to capture the optimal pricing of the TNC platform.  Through a realistic case study, we show that the proposed business model will offer a Pareto improvement that benefits all stakeholders, which leads to higher passenger surplus, higher drivers surplus, higher garage operator surplus, higher platform profit, and reduced traffic congestion.

\end{abstract}

\section{Introduction}
The proliferation of smartphones has enabled transportation network companies (TNCs), such as Uber, Lyft and Didi, to offer real-time matching services for passengers and drivers, which significantly reduces the search frictions of the ride-hailing market. These emerging businesses have deeply transformed the landscape of urban transportation. Today, Uber has completed more than 10 billion trips with 3.9 million active drivers in over 65 countries \cite{Uber_S1}. It was estimated in  \cite{cohen2016using} that Uber generated \$6.8B consumer surplus in the US in 2016. 

Recently, the rapid growth of TNCs starts to exhibit externalities that negatively impact the transportation system. Among the various concerns, a major one is the ``cruising congestion" caused by TNC vehicles. The TNC business model crucially relies on a very short passenger waiting time, which in turn depends on a large number of available but {\em idle} drivers. In New York City, TNC drivers spend more than 40\% of their time empty and cruising for passengers \cite{schaller2017empty}. This underutilziation of vehicles not only  leads to low driver incomes, but also created more cruising traffic that congested the city streets which are already saturated. In Manhattan, for-hire vehicles make up nearly 30\% of all traffic \cite{NYC2019improving}. 

One way to address the cruising congestion is by regulatory intervention \cite{li2019regulating}, \cite{li2020congestion}. In 2019, the New York City Taxi and Limousine Commission (NYCTLC) required TNC platforms to cap their time cruising without passenger below 31\% out of all driving time in Manhattan core at peak hours (currently 41\% industry-wide) \cite{NYC2019improving}. This regulation aims to motivate TNCs to  better utilize their drivers' resources and reduce the cruising congestion. However, it was overturned by the Supreme Court of State of New York \cite{NYC2019crusing}, claiming that the city’s cruising cap was ``arbitrary and capricious``. 

Another way to curb cruising congestion is by parking the vehicles when there is no passenger. Xu et al. \cite{xu2017optimal} studied the allocation of road space to on-street parking for vacant TNC vehicles to reduce cruising congestion. The optimal parking provision strategy was proposed to address the trade-off between reduced cruising traffic and reduced road spaces. Ruch et al. \cite{ruch2019many} considered a mobility-on-demand system with a fleet of vehicles and a number of free parking spaces. A coordinated parking operating policy was proposed to meet the parking capacity constraint with minimum increase in vehicle distance traveled and minimum impact on the service level. Kondor et al. \cite{kondor2018large} investigated the minimum number of parking spaces needed for the mobility-on-demand services and analyzed the relationship between parking provision and traffic congestion. Lam et al \cite{lam2017coordinated} proposed a coordinated parking management strategy for electric  autonomous vehicles to provide vehicle-to-grid services. Jian et al \cite{jian2019unlock} considered an operator that integrates car-sharing and parking sharing platforms to provide bundled services (car and parking) to travelers. 

In this paper, we propose a novel business model that integrates ride-hailing services and parking services into a single ride-sourcing platform (i.e., TNC).  The platform (a) provides interfaces that connect passengers, drivers and garage operators (commercial or private garages); (b) determines the ride fare, driver payment, and parking rates; and (c) executes two matching processes for the market participants, where passengers are matched to TNC vehicles for ride-hailing services, and vacant TNC vehicles are matched to unoccupied parking garages to reduce the cruising cost. Our key observation is as follows: regular parking demands are typically much longer than TNC parking demands (e.g., 1h vs 8min \cite{NYC2019improving}). Based on the SF-Park data \cite{SFpark}, there are inevitably gaps between two long-term parking customers subsequently taking the same parking space. The proposed business model enables the garage owners to fill these gaps using short-term TNC parking demands in real-time, which only takes the spaces that would have been {\em unused} if they are not provided for TNC parking. By unlocking  the potential of these underutilzied parking garages,  the proposed business model will benefit passengers, drivers, garage operators, and the TNC platform without affecting  regular parking demand. 

The key contributions of the paper are summarized below:
\begin{itemize}
    \item We proposed a novel business model that integrated parking services into the TNC platform. Based on the principle of sharing economy, the proposed solution will increase passenger surplus, driver surplus, garage profit, and TNC profit. In the meanwhile, it will reduce the cruising congestion without affecting the regular parking demand. 
    \item We developed a market equilibrium model that captures the incentives of the passengers, drivers, garage operators, and the TNC platform. A  queuing theoretic model is proposed to capture the dynamics of passengers, drivers, and the parking garages, and an optimization model is formulated for the TNC platform to jointly derive the optimal prices of the ride-hailing service and the parking service. 
    \item We validated the proposed model through numerical simulations. The impacts of the proposed parking services on the ride-hailing market are assessed quantitively based on real San Francisco data. 
\end{itemize}

\section{The Business Model}
We consider a transportation market with the TNC platform, a group of passengers and drivers, and a number of  (private or commercial) parking garages distributed across the transportation network.  Based on the principle of sharing economy, the TNC platform offers app-based interfaces that provide {\em on-demand} services to both passengers, drivers and the garage operators. Passengers can request on-demand ride services through the user app. Drivers can log on/off the app to provide ride services depending on their own work schedules. Garage operators can offer or withdraw any number of parking spots at any time depending on their operational strategies.  Each ride is initiated by a service request from the passenger. The platform then matches the passenger with a nearby driver, who delivers the passenger to a pre-specified destination. After the passenger alights, instead of cruising on the street, the driver is dispatched to a nearby parking garage\footnote{We will show that on average, drivers benefit from following the platform dispatch. The spatial heterogeneity among drivers is left for further research. },  where he stays for a few minutes (e.g. 5-8 min \cite{NYC2019improving}) until the next passenger arrives. Throughout this process, a few transactions take place: (a) each passenger pays a ride fare to the platform based on the trip distance and trip time; (b) each driver receives a payment from the platform based on the passenger fare and the commission rate;  and (c) each garage operator receives a payment from the platform based on the accumulated TNC parking time at a per minute rate.  We assume that all rates are determined by the platform, and all transactions are centrally processed  through the cloud. In this case, there is no direct payment transfer among passengers, drivers and garage operators. The system diagram is shown in Figure 1. 
\begin{figure}[bt]%
\centering
\includegraphics[width = 0.8\linewidth]{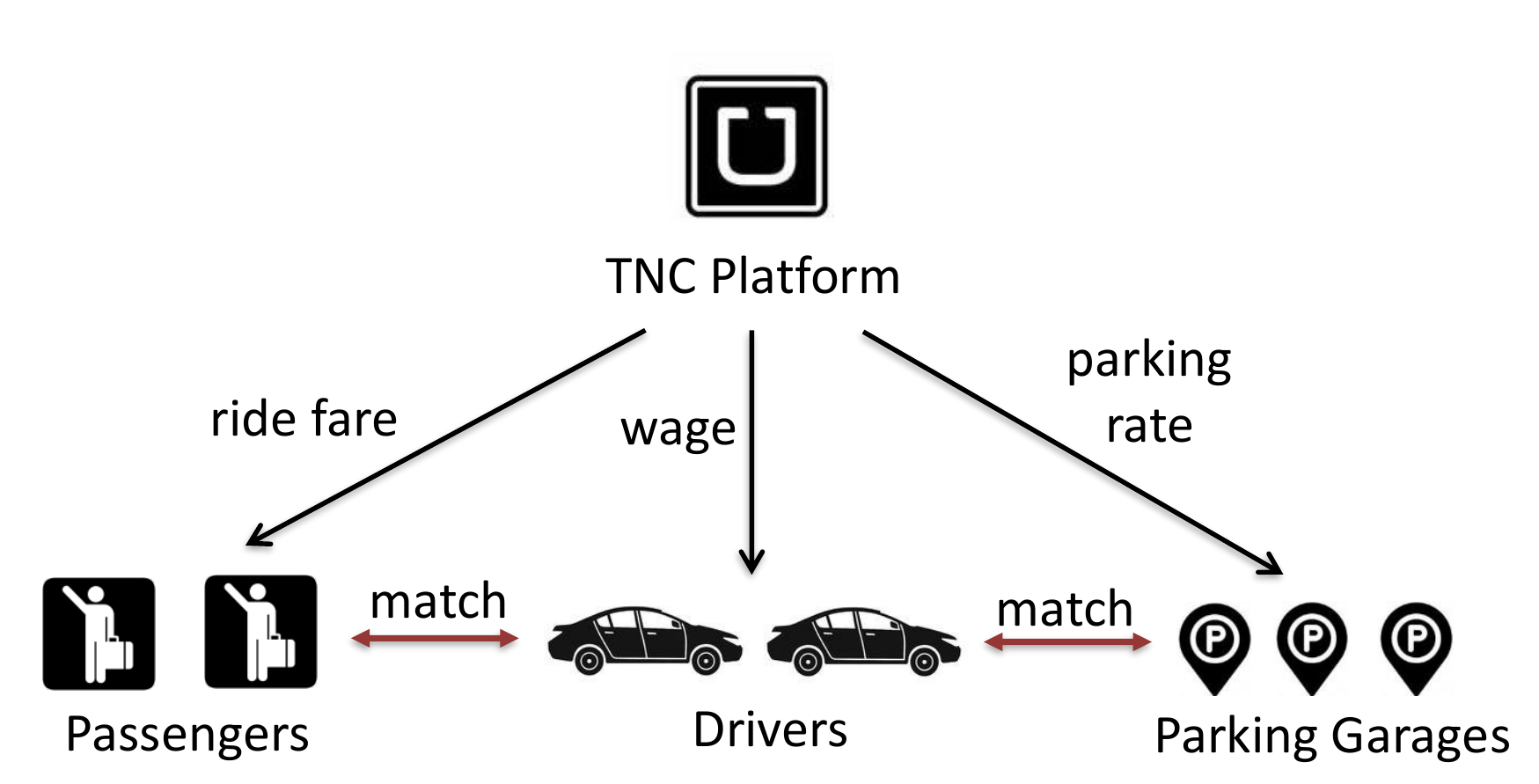}
\caption{The system diagram for the transportation market.}
\label{figdiagram}
\end{figure}
A few important remarks are in order:
\begin{itemize}
\item We argue that the parking service should be provided by the TNC platform rather than other third-party platforms. This is because the sharing-economy model for parking services requires the platform to (a) acquire the real-time location of TNC vehicles, and (b) accumulate a large number of users on both demand side (TNC vehicles) and supply side (parking garage) of the market. TNC is a natural candidate for running the proposed business model\footnote{We will show that TNC makes more profit by offering parking services.}.
\item The proposed business model aims to unlock the potential of the {\em underutilized} parking garages, so that parking services complement the existing garage operation instead of substituting it. It is crucial to observe that each TNC vehicle only parks for 5-8 minutes before the next ride request. This enables the garage operator to fill the gap between regular long-term parking demands (e.g., 2 hours) using the short-term TNC parking demands (e.g., 5-8 min). The garage operator could offer more parking spaces to the TNC platform when it predicts that the regular long-term parking demand is low, and withdraw the parking services from the platform at any time when the regular demand is high. In this case, the parking spaces offered to the TNC platform would have been {\em unoccupied} if it is not reserved for TNC parking. Therefore, the proposed business model will {\em not} impact the regular parking demand\footnote{Garage operators needs to reserve easy-access parking spaces for TNC vehicles to accommodate their short stay.}.
\item  Based on the above discussion, the marginal cost of the parking spaces offered to the TNC platform are minimum. This enables TNC vehicles to enjoy a heavy parking discount compared to other parking demands. The rates of the TNC parking services are determined by the platform, while the rates for other regular parking demands are determined by the garage operator. These two parking services are operated in parallel. The garage operator determines the allocation of the parking spaces in each service to maximize his profit. 
\item The proposed business model requires minimum hardware and software modification on the existing parking infrastructure. The TNC platform will develop an app-based interface for the garage operator. The garage operator will communicate with the TNC platform to enable smooth entry and exit of TNC vehicles. Parking data should be submitted periodically through the interface for financial settlements. 
\end{itemize}

\section{Mathematical Formulation}
Ride-hailing market is typically formulated as a two-sided market where the platform matches the supply side (driver) and the demand side (passenger) to make a profit. When parking service is integrated, the TNC platform faces a much more complicated market structure that involves  two matching processes (i.e., matching passengers to drivers, matching drivers to parking garages) and  three sides (i.e., passenger, drivers, and parking garages). In this case, the platform determines the ride fare, driver payment, and the parking rate to attract both passengers, drivers and the garage operators. These platform decisions will affect the decision making of passengers, drivers and the garage operators, which in turn determines the profit of the TNC platform. Below we present a mathematical model to capture the incentives of these market participants. 

\subsection{The Matching Process}
Consider a ride-hailing market that consists of $N$ drivers, $K$ parking slots, and a group of randomly arriving passengers. Upon arrival, each  passengers immediately joins a queue and waits for he next available driver to be dispatched. After passenger alights, each driver is dispatched to a parking garage (if available)  until the next ride request arrives (driver relocation is neglected in the aggregate model). This can be captured by a continuous-time queuing process where  passengers  are modeled as ``jobs" and drivers are modeled as ``servers". We assume that the arrival process of passengers is Poisson with rate $\lambda>0$, and the service rate is exponential with rate $\mu$. This leads to an M/M/N queue, where the average number of idle TNC vehicles (vehicles without a passenger) is $N_I=N-\lambda/\mu$.

\subsection{Passenger Incentives}
Passengers choose whether to take TNC based on the travel cost of the TNC ride. We define the total travel cost of the TNC ride as the weighted sum of the ride fare  $p_f$ and the waiting time $t_w$:
\begin{equation}
c=\alpha  t_w+p_f.
\end{equation}
where $\alpha$ represents the trade off between time and money, the ride fare $p_f$ is determined by the platform, and the passenger waiting time can be decomposed into two parts: (a) from the ride being requested to a vehicle being dispatched (denoted as $t_1$), and (b) from the vehicle being dispatched to the passenger being picked up (denoted as $t_2$). The former represents the waiting time in the M/M/N queue, and the latter represents the travel time between the passenger and the nearest idle TNC vehicle. We assume that  the total waiting time $t_w=t_1+t_2$  depends on the density of idle TNC vehicles $N_I$. With slight abuse of notation, we denote  $t_w$ as a function $t_w(N-\lambda/\mu)$, then the arrival rate of passengers is captured by:
\begin{equation}
\lambda=F_p\bigg( \alpha t_w(N-\lambda/\mu)+p_f\bigg),
\end{equation}
We assume that $F_p(\cdot)$ is a decreasing and continuously differentiable function, so that higher travel cost $c$ indicates lower arrival rate of TNC passengers.

\subsection{Driver Incentives}
Drivers decide whether to join TNC based on the net hourly wage offered by the platform. The net hourly wage of the TNC driver is defined as the gross wage $w_g$ minus the vehicle expenses. The gross wag $w_g$ is set by the platform \cite{li2019regulating},  determined as the passenger ride fare minus the platform commission. The vehicle expenses consist of a one-time fixed cost (such as vehicle registration, maintenance, insurance) and the variable costs (such as gas expenses, work time, etc). Here we neglect the one-time fixed expenses since the fixed costs are discounted over the long term (e.g. 1 year) and does not affect the short-term driver supply in the TNC market (e.g., 1 hour). 

After the passenger alights, drivers can either cruise at the cost of $l$ (per hour) or park in the garage at the cost of $p_g$ (per hour) subject to parking availability. Let $r$ denote the average utilization rate of the $K$ parking garages, then the net hourly wage of the drivers can be written as:
\begin{equation}
\label{net_wage}
w_n=w_g+  \dfrac{(l-p_g)Kr}{N}
\end{equation}
To derive equation (\ref{net_wage}), we note that $\dfrac{\lambda}{N \mu}$ is the occupancy rate of the TNC vehicles. If no parking service is provided, then the variable cost for all vehicles is $l$ per hour. However, the TNC vehicle would save $l-p_g$ (per hour parked) if the parking rate is lower than the cruising cost. Since the occupancy rate of the parking garages are $r$, the total saving for all TNC drivers are $(l-p_g)Kr$. Therefore each drivers receives a saving of  $\dfrac{(l-p_g)Kr}{N}$. 

The total number of drivers $N$ is a function on the net hourly wage:
\begin{equation}
\label{driver_cost}
N=F_d\left( w_g+ \dfrac{(l-p_g)Kr}{N}\right). 
\end{equation}
We assume that $F_d(\cdot)$ is an increasing function, so that higher net wage indicates more supply of TNC drivers. Without loss of generality,  the vehicle variable cost $-l$ is neglected in  (\ref{driver_cost}) as it is a constant.

\subsection{Garage Operator Incentives}
The garage operator decides the number of parking spaces offered to the TNC platform based on the TNC parking rates and the statistics of the regular parking demand. In practice, the garage operator should predict the regular parking demand in the next a few minutes, evaluate how many parking spaces will be unoccupied, and offer these parking spaces to the TNC platform. To reduce the exposure, we assume that this  underlying decision process is captured by a function that specifies the relation between the TNC parking rate $p_g$ and the parking supply $K$:
\begin{equation}
K=F_g(rp_g)
\end{equation}
where $rp_g$ is the per hour earning of a parking space offered to the TNC platform. We assume that $F_g$ is an increasing function so that higher earnings attract more parking supply. 

It is important to note that the utilization rate of TNC parking space, $r$,  depends on the passenger arrival rate $\lambda$, the number of drivers $N$, and the number of parking spaces $K$. With slight abuse of notation, denote $r$ as $r(\lambda, N, K)$. We note that when there are $i$ idle drivers in the system and $i\geq K$, all parking space will be occupied, and $r=1$. However, when $i<K$, only $i$ out of $K$ parking spaces are occupied, and the average utilization rate satisfies $r=i/K$. In addition, we note that the passenger-driver matching process is an M/M/N queue. Therefore, the probability that there are $i$ {\em idle} drivers, denoted as $X_i$, can be derived as:
\begin{align}
 \begin{cases}
& X_0=\pi_0 \dfrac{(N\rho)^N}{N!(1-\rho)} \\
& X_i=\pi_0 \dfrac{(N\rho)^{N-i}}{(N-i)!}, \quad \forall 1\leq i\leq N-1\\
& X_N=\pi_0,
\end{cases}   
\end{align}
where $\rho=\dfrac{\lambda}{N \mu}$ is the occupancy rate of TNC vehicles, and $\pi_0$ is the probability of all drivers being idle in the ride-hailing market. 
The expected utilization rate of parking spaces on the TNC platform is 
\begin{equation}
r(\lambda, N, K)=\sum_{i=1}^N X_i \min \left(1,\dfrac{i}{K}\right),
\end{equation}
where $\min(a,b)$ represents the min of $a$ and $b$.

\subsection{Profit Maximization of the TNC Platform}
Consider a TNC platform that charges commission to make a profit. In this paper, we assume that the platform charges a commission on the ride-hailing services, but does not charge a commission on  the parking services. We argue that the TNC platform will refrain from taking commission from the parking services at least in the early stage of the proposed business model. In the next section, we will show that the TNC platform substantially benefits from the proposed business model even if no commission charge is directly collected from  the parking sector. 

In each hour, the platform receives $\lambda p_f$ from passengers,  pays $w_gN$ to drivers, and keeps the difference to make a profit. The platform's profit maximization problem can be cast as:
\begin{align}
\label{TNCprofit_max}
&\max_{p_f, w_g, p_g} \lambda p_f-w_gN \\
& 
\begin{cases}
\label{marekt_constraints}
\lambda=F_p\bigg( \alpha t_w(N-\lambda/\mu)+p_f\bigg), \\
N=F_d\left( w_g +  \dfrac{(l-p_g)Kr(\lambda, N, K)}{N}\right)  \\
K=F_g\left(p_gr(\lambda, N, K)\right).
\end{cases}
\end{align}
The overall problem (\ref{TNCprofit_max}) is non-convex due to the nonlinear constraints. However,  since the dimension of the problem is small, we can solve the globally optimal solution to (\ref{TNCprofit_max}) via enumerations within sub-second. 

\section{Main Results}
This section discusses the feasibility  of the proposed business model. We first demonstrate the benefits of the proposed model  through an numerical example using real data from San Francisco, then we summarize the theoretical result in Theorem 1, which shows that integrating parking services to the TNC platform will offer a Pareto improvement, leading to higher passenger surplus, higher driver surplus, higher platform profit, and less traffic congestion. 

\subsection{Numerical Example}
Passengers choose their transport mode based on the total travel cost $c$. This can be captured by a logit model, which can be written as:
\begin{equation}
\lambda=\lambda_0\dfrac{e^{-\epsilon c}}{e^{-\epsilon c}+e^{-\epsilon c_0}},
\end{equation}
where $\lambda_0$ is the number of potential passengers (all passengers regardless of their travel mode), and $\epsilon$ and $c_0$ are the model parameters. In our previous work \cite{li2019regulating}, we proved that the passenger waiting time $t_w$ is inversely proportional to the square root of the number of idle vehicles:
\begin{equation}
t_w=\dfrac{M}{\sqrt{N-\lambda/\mu}}.
\end{equation}
The underlying intuition of this results is very simple: assume that all idle vehicles are uniformly distributed across the city. The passenger waiting time depends on the distance between the passenger and the closest idle vehicle, which is further proportional to the distance between any two nearby idle vehicles. Intuitively, the average distance between any two nearby idle vehicles is inversely proportional to the square root of the number of idle vehicles. Our result is consistent with the empirical data reported in \cite{korolko2018dynamic}.

\begin{figure*}[bt]
\begin{minipage}[b]{0.32\linewidth}
\centering
%

%
%
\definecolor{mycolor1}{rgb}{0.00000,0.44700,0.74100}%
\begin{tikzpicture}

\begin{axis}[%
width=1.7in,
height=0.8in,
at={(1.358in,0.0in)},
scale only axis,
xmin=0,
xmax=3000,
xtick={0,1000,2000,3000},
xticklabels={{0},{1e3},{2e3},{3e3}},
xlabel style={font=\color{white!15!black}},
xlabel={Parking Supply},
ymin=0.4,
ymax=1.1,
ylabel style={font=\color{white!15!black}},
ylabel={Garage Occupancy},
axis background/.style={fill=white},
legend style={legend cell align=left, align=left, draw=white!15!black}
]
\addplot [color=black, line width=1.0pt]
  table[row sep=crcr]{%
0	1\\
30.3030303030303	0.999999999999999\\
60.6060606060606	0.999999999999999\\
90.9090909090909	1\\
121.212121212121	0.999999999999999\\
151.515151515152	1\\
181.818181818182	1\\
212.121212121212	0.999999999999999\\
242.424242424242	0.999999999999999\\
272.727272727273	1\\
303.030303030303	0.999999999999999\\
333.333333333333	0.999999999999999\\
363.636363636364	1\\
393.939393939394	0.999999999999999\\
424.242424242424	0.999999999999999\\
454.545454545455	0.999999999999999\\
484.848484848485	0.999999999999999\\
515.151515151515	0.999999999999999\\
545.454545454545	0.999999999999999\\
575.757575757576	0.999999999999999\\
606.060606060606	1\\
636.363636363636	0.999999999999999\\
666.666666666667	0.999999999999999\\
696.969696969697	1\\
727.272727272727	0.999999999999999\\
757.575757575758	1\\
787.878787878788	0.999999999999999\\
818.181818181818	0.999999999999999\\
848.484848484848	1\\
878.787878787879	0.999999999999999\\
909.090909090909	0.999999999999999\\
939.393939393939	0.999999999999999\\
969.69696969697	1\\
1000	1\\
1030.30303030303	1\\
1060.60606060606	1\\
1090.90909090909	0.999999999999957\\
1121.21212121212	0.999999999995553\\
1151.51515151515	0.999999999666513\\
1181.81818181818	0.999999986151913\\
1212.12121212121	0.999999632035739\\
1242.42424242424	0.999994138832546\\
1272.72727272727	0.999933383647717\\
1303.0303030303	0.999617582296976\\
1333.33333333333	0.998780834886491\\
1363.63636363636	0.997262030425345\\
1393.93939393939	0.995205816565024\\
1424.24242424242	0.992603594939557\\
1454.54545454545	0.989211644214574\\
1484.84848484848	0.985429611240903\\
1515.15151515152	0.980640570787269\\
1545.45454545455	0.975394828522893\\
1575.75757575758	0.968088814122067\\
1606.06060606061	0.958518578838384\\
1636.36363636364	0.946512883788967\\
1666.66666666667	0.931462247775615\\
1696.9696969697	0.915142671952558\\
1727.27272727273	0.899264907080643\\
1757.57575757576	0.88339131971815\\
1787.87878787879	0.868576542884619\\
1818.18181818182	0.85425042214588\\
1848.48484848485	0.840389216571247\\
1878.78787878788	0.826509587352451\\
1909.09090909091	0.813527165947676\\
1939.39393939394	0.800946319594986\\
1969.69696969697	0.788329168400518\\
2000	0.776509923322471\\
2030.30303030303	0.765039876689795\\
2060.60606060606	0.753520493756236\\
2090.90909090909	0.742714811886699\\
2121.21212121212	0.732214647553755\\
2151.51515151515	0.721655758358655\\
2181.81818181818	0.711738622885805\\
2212.12121212121	0.702090381284092\\
2242.42424242424	0.692700218472284\\
2272.72727272727	0.683242829311837\\
2303.0303030303	0.674346878479282\\
2333.33333333333	0.665679572411582\\
2363.63636363636	0.656940973667097\\
2393.93939393939	0.64871260540799\\
2424.24242424242	0.640687816609537\\
2454.54545454545	0.632589059695589\\
2484.84848484848	0.62495588167512\\
2515.15151515152	0.617504698976039\\
2545.45454545455	0.61022911621593\\
2575.75757575758	0.602877669522841\\
2606.06060606061	0.59594075334621\\
2636.36363636364	0.58916164449167\\
2666.66666666667	0.582306191548582\\
2696.9696969697	0.575832053467251\\
2727.27272727273	0.569500281430698\\
2757.57575757576	0.563092263143226\\
2787.87878787879	0.557036108915024\\
2818.18181818182	0.551108820471394\\
2848.48484848485	0.545306348597319\\
2878.78787878788	0.539428415362342\\
2909.09090909091	0.533868053223433\\
2939.39393939394	0.528421168110818\\
2969.69696969697	0.522899768714387\\
3000	0.517673282214982\\
};

\end{axis}

\end{tikzpicture} 
\vspace*{-0.3in}
\caption{Average occupancy rate of the parking garages as a function of $K$. }
\label{figure1}
\end{minipage}
\begin{minipage}[b]{0.005\linewidth}
\hfill
\end{minipage}
\begin{minipage}[b]{0.32\linewidth}
\centering
%
%
\definecolor{mycolor1}{rgb}{0.00000,0.44700,0.74100}%
\begin{tikzpicture}

\begin{axis}[%
width=1.7in,
height=0.8in,
at={(1.358in,0.0in)},
scale only axis,
xmin=0,
xmax=3000,
xtick={0,1000,2000,3000},
xticklabels={{0},{1e3},{2e3},{3e3}},
xlabel style={font=\color{white!15!black}},
xlabel={Parking Supply},
ymin=0,
ymax=1,
ylabel style={font=\color{white!15!black}},
ylabel={Parking Ratio},
axis background/.style={fill=white},
legend style={legend cell align=left, align=left, draw=white!15!black}
]
\addplot [color=black, line width=1.0pt]
  table[row sep=crcr]{%
0	0\\
30.3030303030303	0.0218441031391733\\
60.6060606060606	0.0444163430496525\\
90.9090909090909	0.0662604487973922\\
121.212121212121	0.0881045493279991\\
151.515151515152	0.110676793595644\\
181.818181818182	0.132520897594784\\
212.121212121212	0.154364995516825\\
242.424242424242	0.176209098655998\\
272.727272727273	0.198781346392177\\
303.030303030303	0.22062544170565\\
333.333333333333	0.242469544844824\\
363.636363636364	0.265041795189569\\
393.939393939394	0.286885887894476\\
424.242424242424	0.30872999103365\\
454.545454545455	0.331302230944128\\
484.848484848485	0.353146334083302\\
515.151515151515	0.374990437222475\\
545.454545454545	0.396834540361648\\
575.757575757576	0.419406780272128\\
606.060606060606	0.441250900782635\\
636.363636363636	0.463094986550474\\
666.666666666667	0.485667226460953\\
696.969696969697	0.507511349580026\\
727.272727272727	0.5293554327393\\
757.575757575758	0.551927694378278\\
787.878787878788	0.573771775788952\\
818.181818181818	0.595615878928125\\
848.484848484848	0.617460006375699\\
878.787878787879	0.640032221977778\\
909.090909090909	0.661876325116951\\
939.393939393939	0.683720428256124\\
969.69696969697	0.706292695972204\\
1000	0.728136799971345\\
1030.30303030303	0.749980903970483\\
1060.60606060606	0.772553144769596\\
1090.90909090909	0.794397248768693\\
1121.21212121212	0.816241320629277\\
1151.51515151515	0.838813560206426\\
1181.81818181818	0.860657271281824\\
1212.12121212121	0.882491966352795\\
1242.42424242424	0.903799812558336\\
1272.72727272727	0.925788786330662\\
1303.0303030303	0.944588171096499\\
1333.33333333333	0.958944515226699\\
1363.63636363636	0.969849117068624\\
1393.93939393939	0.977634254367522\\
1424.24242424242	0.983557465754129\\
1454.54545454545	0.988421353085957\\
1484.84848484848	0.991933708301207\\
1515.15151515152	0.994795775280319\\
1545.45454545455	0.99677310251996\\
1575.75757575758	0.998364033103444\\
1606.06060606061	0.999366782308765\\
1636.36363636364	0.999834171617843\\
1666.66666666667	0.99997740578157\\
1696.9696969697	0.999998407499103\\
1727.27272727273	0.999999932912118\\
1757.57575757576	0.999999998542185\\
1787.87878787879	0.999999999979257\\
1818.18181818182	0.99999999999983\\
1848.48484848485	0.999999999999998\\
1878.78787878788	1\\
1909.09090909091	1\\
1939.39393939394	0.999999999999999\\
1969.69696969697	1\\
2000	1\\
2030.30303030303	1\\
2060.60606060606	1\\
2090.90909090909	1\\
2121.21212121212	1\\
2151.51515151515	1\\
2181.81818181818	1\\
2212.12121212121	1\\
2242.42424242424	1\\
2272.72727272727	1\\
2303.0303030303	0.999999999999998\\
2333.33333333333	1\\
2363.63636363636	1\\
2393.93939393939	0.999999999999999\\
2424.24242424242	1\\
2454.54545454545	1\\
2484.84848484848	1\\
2515.15151515152	1\\
2545.45454545455	1\\
2575.75757575758	1\\
2606.06060606061	0.999999999999999\\
2636.36363636364	1\\
2666.66666666667	1\\
2696.9696969697	1\\
2727.27272727273	1\\
2757.57575757576	1\\
2787.87878787879	1\\
2818.18181818182	1\\
2848.48484848485	0.999999999999999\\
2878.78787878788	1\\
2909.09090909091	1\\
2939.39393939394	1\\
2969.69696969697	0.999999999999999\\
3000	1\\
};

\end{axis}

\end{tikzpicture}
\vspace*{-0.3in}
\caption{Ratio of parked TNC vehicles out of all idle TNC vehicles as a function of $K$.} 
\label{figure2}
\end{minipage}
\begin{minipage}[b]{0.005\linewidth}
\hfill 
\end{minipage}
\begin{minipage}[b]{0.32\linewidth}
\centering
%
%
\definecolor{mycolor1}{rgb}{0.00000,0.44700,0.74100}%
\begin{tikzpicture}

\begin{axis}[%
width=1.7in,
height=0.8in,
at={(1.358in,0.0in)},
scale only axis,
xmin=0,
xmax=3000,
xtick={0,1000,2000,3000},
xticklabels={{0},{1e3},{2e3},{3e3}},
xlabel style={font=\color{white!15!black}},
xlabel={Parking Supply},
ymin=150,
ymax=160,
ylabel style={font=\color{white!15!black}},
ylabel={Passenger Arrivals},
axis background/.style={fill=white},
legend style={legend cell align=left, align=left, draw=white!15!black}
]
\addplot [color=black, line width=1.0pt]
  table[row sep=crcr]{%
0	150.626318463016\\
30.3030303030303	150.626318463016\\
60.6060606060606	150.626318463016\\
90.9090909090909	150.626323307336\\
121.212121212121	150.626318463016\\
151.515151515152	150.626323307336\\
181.818181818182	150.626323307336\\
212.121212121212	150.626318463016\\
242.424242424242	150.626318463016\\
272.727272727273	150.626323307336\\
303.030303030303	150.626318463016\\
333.333333333333	150.626318463016\\
363.636363636364	150.626323307336\\
393.939393939394	150.626318463016\\
424.242424242424	150.626318463016\\
454.545454545455	150.626318463016\\
484.848484848485	150.626318463016\\
515.151515151515	150.626318463016\\
545.454545454545	150.626318463016\\
575.757575757576	150.626318463016\\
606.060606060606	150.626323307336\\
636.363636363636	150.626318463016\\
666.666666666667	150.626318463016\\
696.969696969697	150.626323307336\\
727.272727272727	150.626318463016\\
757.575757575758	150.626323307336\\
787.878787878788	150.626318463016\\
818.181818181818	150.626318463016\\
848.484848484848	150.626323307336\\
878.787878787879	150.626318463016\\
909.090909090909	150.626318463016\\
939.393939393939	150.626318463016\\
969.69696969697	150.626323307336\\
1000	150.626323307336\\
1030.30303030303	150.626323307336\\
1060.60606060606	150.626323307336\\
1090.90909090909	150.626323307336\\
1121.21212121212	150.626318463016\\
1151.51515151515	150.626318463016\\
1181.81818181818	150.626264368112\\
1212.12121212121	150.625004042057\\
1242.42424242424	150.642464207532\\
1272.72727272727	150.664615414284\\
1303.0303030303	150.851175089661\\
1333.33333333333	151.35476580727\\
1363.63636363636	152.066886381479\\
1393.93939393939	152.84090823793\\
1424.24242424242	153.746404480561\\
1454.54545454545	154.560756202228\\
1484.84848484848	155.460572666508\\
1515.15151515152	156.173579053206\\
1545.45454545455	157.132189303396\\
1575.75757575758	157.828064591002\\
1606.06060606061	158.270213248281\\
1636.36363636364	158.687229027215\\
1666.66666666667	158.955124794366\\
1696.9696969697	158.935911086754\\
1727.27272727273	158.933150248421\\
1757.57575757576	158.935723095633\\
1787.87878787879	158.93456396147\\
1818.18181818182	158.933450356018\\
1848.48484848485	158.932373005569\\
1878.78787878788	158.934864071195\\
1909.09090909091	158.933801044118\\
1939.39393939394	158.932764155606\\
1969.69696969697	158.935138048165\\
2000	158.934116326736\\
2030.30303030303	158.933119900487\\
2060.60606060606	158.935382520241\\
2090.90909090909	158.934401261723\\
2121.21212121212	158.933450356018\\
2151.51515151515	158.935610975452\\
2181.81818181818	158.934670180029\\
2212.12121212121	158.933750464066\\
2242.42424242424	158.932856042256\\
2272.72727272727	158.93491465153\\
2303.0303030303	158.934025282506\\
2333.33333333333	158.93316542239\\
2363.63636363636	158.935138048165\\
2393.93939393939	158.934284084268\\
2424.24242424242	158.933450356018\\
2454.54545454545	158.935346270914\\
2484.84848484848	158.934518439247\\
2515.15151515152	158.933715058036\\
2545.45454545455	158.932926854016\\
2575.75757575758	158.934741835438\\
2606.06060606061	158.933958685365\\
2636.36363636364	158.933196613329\\
2666.66666666667	158.934944999738\\
2696.9696969697	158.934187981946\\
2727.27272727273	158.933450356018\\
2757.57575757576	158.935138048165\\
2787.87878787879	158.934401261723\\
2818.18181818182	158.933684710016\\
2848.48484848485	158.932982491844\\
2878.78787878788	158.934604425676\\
2909.09090909091	158.933908105272\\
2939.39393939394	158.933221903284\\
2969.69696969697	158.934792415742\\
3000	158.934116326736\\
};

\end{axis}

\end{tikzpicture}%
\vspace*{-0.3in}
\caption{Passenger arrival rate (per minute) as a function of $K$.}
\label{figure3}
\end{minipage}
\end{figure*}

\begin{figure*}[bt]
\begin{minipage}[b]{0.32\linewidth}
\centering
%
%
\definecolor{mycolor1}{rgb}{0.00000,0.44700,0.74100}%
\definecolor{mycolor2}{rgb}{0.85000,0.32500,0.09800}%
\begin{tikzpicture}

\begin{axis}[%
width=1.7in,
height=0.8in,
at={(1.358in,0.0in)},
scale only axis,
xmin=0,
xmax=3000,
xtick={0,1000,2000,3000},
xticklabels={{0},{1e3},{2e3},{3e3}},
xlabel style={font=\color{white!15!black}},
xlabel={Parking Supply},
ymin=1500,
ymax=3500,
ytick={1000,2000,3000},
yticklabels={{1e3},{2e3},{3e3}},
ylabel style={font=\color{white!15!black}},
ylabel={Number of Driver},
axis background/.style={fill=white},
legend style={at={(0.6,0.35)}, anchor=south west, legend cell align=left, align=left, draw=white!12!black}
]
\addplot [color=black, line width=1.0pt]
  table[row sep=crcr]{%
0	3053.13490710308\\
30.3030303030303	3053.13490710308\\
60.6060606060606	3053.13490710308\\
90.9090909090909	3053.13497304011\\
121.212121212121	3053.13490710308\\
151.515151515152	3053.13497304011\\
181.818181818182	3053.13497304011\\
212.121212121212	3053.13490710308\\
242.424242424242	3053.13490710308\\
272.727272727273	3053.13497304011\\
303.030303030303	3053.13490710308\\
333.333333333333	3053.13490710308\\
363.636363636364	3053.13497304011\\
393.939393939394	3053.13490710308\\
424.242424242424	3053.13490710308\\
454.545454545455	3053.13490710308\\
484.848484848485	3053.13490710308\\
515.151515151515	3053.13490710308\\
545.454545454545	3053.13490710308\\
575.757575757576	3053.13490710308\\
606.060606060606	3053.13497304011\\
636.363636363636	3053.13490710308\\
666.666666666667	3053.13490710308\\
696.969696969697	3053.13497304011\\
727.272727272727	3053.13490710308\\
757.575757575758	3053.13497304011\\
787.878787878788	3053.13490710308\\
818.181818181818	3053.13490710308\\
848.484848484848	3053.13497304011\\
878.787878787879	3053.13490710308\\
909.090909090909	3053.13490710308\\
939.393939393939	3053.13490710308\\
969.69696969697	3053.13497304011\\
1000	3053.13497304011\\
1030.30303030303	3053.13497304011\\
1060.60606060606	3053.13497304011\\
1090.90909090909	3053.13497304011\\
1121.21212121212	3053.13490710308\\
1151.51515151515	3053.13490710308\\
1181.81818181818	3053.13409387976\\
1212.12121212121	3053.11568649251\\
1242.42424242424	3053.50000448542\\
1272.72727272727	3054.5000470331\\
1303.0303030303	3060.50002310462\\
1333.33333333333	3075.50002637136\\
1363.63636363636	3097.50000052376\\
1393.93939393939	3122.50001780637\\
1424.24242424242	3150.500054403\\
1454.54545454545	3178.50003451991\\
1484.84848484848	3207.50001701044\\
1515.15151515152	3233.50000993425\\
1545.45454545455	3262.50005899765\\
1575.75757575758	3286.50001058515\\
1606.06060606061	3303.5000516127\\
1636.36363636364	3316.50002455102\\
1666.66666666667	3323.5000455143\\
1696.9696969697	3323.5000455143\\
1727.27272727273	3323.5000455143\\
1757.57575757576	3323.5000455143\\
1787.87878787879	3323.5000455143\\
1818.18181818182	3323.5000455143\\
1848.48484848485	3323.5000455143\\
1878.78787878788	3323.5000455143\\
1909.09090909091	3323.5000455143\\
1939.39393939394	3323.5000455143\\
1969.69696969697	3323.5000455143\\
2000	3323.5000455143\\
2030.30303030303	3323.5000455143\\
2060.60606060606	3323.5000455143\\
2090.90909090909	3323.5000455143\\
2121.21212121212	3323.5000455143\\
2151.51515151515	3323.5000455143\\
2181.81818181818	3323.5000455143\\
2212.12121212121	3323.5000455143\\
2242.42424242424	3323.5000455143\\
2272.72727272727	3323.5000455143\\
2303.0303030303	3323.5000455143\\
2333.33333333333	3323.5000455143\\
2363.63636363636	3323.5000455143\\
2393.93939393939	3323.5000455143\\
2424.24242424242	3323.5000455143\\
2454.54545454545	3323.5000455143\\
2484.84848484848	3323.5000455143\\
2515.15151515152	3323.5000455143\\
2545.45454545455	3323.5000455143\\
2575.75757575758	3323.5000455143\\
2606.06060606061	3323.5000455143\\
2636.36363636364	3323.5000455143\\
2666.66666666667	3323.5000455143\\
2696.9696969697	3323.5000455143\\
2727.27272727273	3323.5000455143\\
2757.57575757576	3323.5000455143\\
2787.87878787879	3323.5000455143\\
2818.18181818182	3323.5000455143\\
2848.48484848485	3323.5000455143\\
2878.78787878788	3323.5000455143\\
2909.09090909091	3323.5000455143\\
2939.39393939394	3323.5000455143\\
2969.69696969697	3323.5000455143\\
3000	3323.5000455143\\
};
\addlegendentry{$N$}

\addplot [color=mycolor1, dashed, line width=1.0pt]
  table[row sep=crcr]{%
0	3053.13490710308\\
30.3030303030303	3023.10522089543\\
60.6060606060606	2992.07454514752\\
90.9090909090909	2962.04492049853\\
121.212121212121	2932.01517273223\\
151.515151515152	2900.98455560802\\
181.818181818182	2870.95486795695\\
212.121212121212	2840.92512456902\\
242.424242424242	2810.89543836137\\
272.727272727273	2779.86481541538\\
303.030303030303	2749.83507640582\\
333.333333333333	2719.80539019817\\
363.636363636364	2688.7747628738\\
393.939393939394	2658.74502824262\\
424.242424242424	2628.71534203497\\
454.545454545455	2597.68466628707\\
484.848484848485	2567.65498007942\\
515.151515151515	2537.62529387177\\
545.454545454545	2507.59560766412\\
575.757575757576	2476.56493191621\\
606.060606060606	2446.53528248851\\
636.363636363636	2416.50555950092\\
666.666666666667	2385.47488375301\\
696.969696969697	2355.44522994693\\
727.272727272727	2325.41551133771\\
757.575757575758	2294.38486505643\\
787.878787878788	2264.35514938216\\
818.181818181818	2234.32546317451\\
848.484848484848	2204.29580210322\\
878.787878787879	2173.26510121896\\
909.090909090909	2143.23541501131\\
939.393939393939	2113.20572880366\\
969.69696969697	2082.1750723222\\
1000	2052.14538467113\\
1030.30303030303	2022.11569702007\\
1060.60606060606	1991.08501978063\\
1090.90909090909	1961.05533212962\\
1121.21212121212	1931.02563248344\\
1151.51515151515	1899.99495719378\\
1181.81818181818	1869.96517787025\\
1212.12121212121	1839.93387337715\\
1242.42424242424	1810.85430297023\\
1272.72727272727	1780.92411177117\\
1303.0303030303	1757.35852873724\\
1333.33333333333	1743.54959548807\\
1363.63636363636	1736.7624982402\\
1393.93939393939	1734.83069459843\\
1424.24242424242	1736.8075005872\\
1454.54545454545	1739.10984944352\\
1484.84848484848	1744.17454826081\\
1515.15151515152	1747.99124623963\\
1545.45454545455	1755.77931141052\\
1575.75757575758	1761.15507978781\\
1606.06060606061	1764.55747931025\\
1636.36363636364	1768.48587070478\\
1666.66666666667	1771.24932039304\\
1696.9696969697	1771.00262446702\\
1727.27272727273	1770.9694926384\\
1757.57575757576	1770.9980596146\\
1787.87878787879	1770.98514131715\\
1818.18181818182	1770.97273253875\\
1848.48484848485	1770.96072777634\\
1878.78787878788	1770.98848536474\\
1909.09090909091	1770.97664020589\\
1939.39393939394	1770.96508630532\\
1969.69696969697	1770.99153825098\\
2000	1770.98015335506\\
2030.30303030303	1770.96905031972\\
2060.60606060606	1770.9942623684\\
2090.90909090909	1770.98332834491\\
2121.21212121212	1770.97273253848\\
2151.51515151515	1770.99680801218\\
2181.81818181818	1770.98632486318\\
2212.12121212121	1770.97607659959\\
2242.42424242424	1770.96611018514\\
2272.72727272727	1770.98904897419\\
2303.0303030303	1770.97913886221\\
2333.33333333333	1770.96955756378\\
2363.63636363636	1770.99153825098\\
2393.93939393939	1770.98202265327\\
2424.24242424242	1770.97273253848\\
2454.54545454545	1770.99385844733\\
2484.84848484848	1770.98463403733\\
2515.15151515152	1770.97568207526\\
2545.45454545455	1770.96689923046\\
2575.75757575758	1770.98712330916\\
2606.06060606061	1770.97839677978\\
2636.36363636364	1770.96990511996\\
2666.66666666667	1770.98938713993\\
2696.9696969697	1770.98095179882\\
2727.27272727273	1770.97273253848\\
2757.57575757576	1770.99153825098\\
2787.87878787879	1770.98332834491\\
2818.18181818182	1770.97534391161\\
2848.48484848485	1770.96751919484\\
2878.78787878788	1770.98559217182\\
2909.09090909091	1770.97783317303\\
2939.39393939394	1770.97018692231\\
2969.69696969697	1770.98768691827\\
3000	1770.98015335506\\
};
\addlegendentry{$N_{on}$}

\end{axis}
\end{tikzpicture}%
\vspace*{-0.3in}
\caption{Total number of drivers and number of on-road drivers ($N_{on}$) as a function of $K$. }
\label{figure4}
\end{minipage}
\begin{minipage}[b]{0.005\linewidth}
\hfill
\end{minipage}
\begin{minipage}[b]{0.32\linewidth}
\centering
%
%
\definecolor{mycolor1}{rgb}{0.00000,0.44700,0.74100}%
\begin{tikzpicture}

\begin{axis}[%
width=1.7in,
height=0.8in,
at={(1.358in,0.0in)},
scale only axis,
xmin=0,
xmax=3000,
xtick={0,1000,2000,3000},
xticklabels={{0},{1e3},{2e3},{3e3}},
xlabel style={font=\color{white!15!black}},
xlabel={Parking Supply},
ymin=28.4,
ymax=29,
ylabel style={font=\color{white!15!black}},
ylabel={Passenger Cost},
axis background/.style={fill=white},
legend style={legend cell align=left, align=left, draw=white!15!black}
]
\addplot [color=black, line width=1.0pt]
  table[row sep=crcr]{%
0	28.8304496365501\\
30.3030303030303	28.8304496365501\\
60.6060606060606	28.8304496365501\\
90.9090909090909	28.8304494416414\\
121.212121212121	28.8304496365501\\
151.515151515152	28.8304494416414\\
181.818181818182	28.8304494416414\\
212.121212121212	28.8304496365501\\
242.424242424242	28.8304496365501\\
272.727272727273	28.8304494416414\\
303.030303030303	28.8304496365501\\
333.333333333333	28.8304496365501\\
363.636363636364	28.8304494416414\\
393.939393939394	28.8304496365501\\
424.242424242424	28.8304496365501\\
454.545454545455	28.8304496365501\\
484.848484848485	28.8304496365501\\
515.151515151515	28.8304496365501\\
545.454545454545	28.8304496365501\\
575.757575757576	28.8304496365501\\
606.060606060606	28.8304494416414\\
636.363636363636	28.8304496365501\\
666.666666666667	28.8304496365501\\
696.969696969697	28.8304494416414\\
727.272727272727	28.8304496365501\\
757.575757575758	28.8304494416414\\
787.878787878788	28.8304496365501\\
818.181818181818	28.8304496365501\\
848.484848484848	28.8304494416414\\
878.787878787879	28.8304496365501\\
909.090909090909	28.8304496365501\\
939.393939393939	28.8304496365501\\
969.69696969697	28.8304494416414\\
1000	28.8304494416414\\
1030.30303030303	28.8304494416414\\
1060.60606060606	28.8304494416414\\
1090.90909090909	28.8304494416414\\
1121.21212121212	28.8304496365501\\
1151.51515151515	28.8304496365501\\
1181.81818181818	28.8304518130303\\
1212.12121212121	28.8305025217711\\
1242.42424242424	28.8298000492425\\
1272.72727272727	28.8289089375991\\
1303.0303030303	28.8214081144813\\
1333.33333333333	28.8011981718159\\
1363.63636363636	28.7727122325674\\
1393.93939393939	28.7418721402971\\
1424.24242424242	28.7059531222314\\
1454.54545454545	28.6737948669394\\
1484.84848484848	28.6384197516134\\
1515.15151515152	28.6105054125901\\
1545.45454545455	28.5731362234351\\
1575.75757575758	28.5461233119277\\
1606.06060606061	28.5290091479956\\
1636.36363636364	28.5129027539642\\
1666.66666666667	28.5025736326899\\
1696.9696969697	28.5033139833738\\
1727.27272727273	28.5034203710267\\
1757.57575757576	28.5033212274797\\
1787.87878787879	28.5033658940515\\
1818.18181818182	28.5034088064452\\
1848.48484848485	28.5034503219934\\
1878.78787878788	28.50335432947\\
1909.09090909091	28.503395292777\\
1939.39393939394	28.5034352490558\\
1969.69696969697	28.5033437719166\\
2000	28.5033831434694\\
2030.30303030303	28.5034215404788\\
2060.60606060606	28.5033343513306\\
2090.90909090909	28.5033721636139\\
2121.21212121212	28.5034088064452\\
2151.51515151515	28.5033255479553\\
2181.81818181818	28.5033618009693\\
2212.12121212121	28.5033972418637\\
2242.42424242424	28.5034317082148\\
2272.72727272727	28.5033523803832\\
2303.0303030303	28.5033866518256\\
2333.33333333333	28.5034197863007\\
2363.63636363636	28.5033437719166\\
2393.93939393939	28.5033766789983\\
2424.24242424242	28.5034088064452\\
2454.54545454545	28.5033357481761\\
2484.84848484848	28.5033676482296\\
2515.15151515152	28.5033986062245\\
2545.45454545455	28.5034289794933\\
2575.75757575758	28.503359039763\\
2606.06060606061	28.5033892181232\\
2636.36363636364	28.5034185843639\\
2666.66666666667	28.5033512109312\\
2696.9696969697	28.5033803822632\\
2727.27272727273	28.5034088064452\\
2757.57575757576	28.5033437719166\\
2787.87878787879	28.5033721636139\\
2818.18181818182	28.5033997756765\\
2848.48484848485	28.5034268354979\\
2878.78787878788	28.5033643347821\\
2909.09090909091	28.50339116721\\
2939.39393939394	28.5034176098205\\
2969.69696969697	28.5033570906762\\
3000	28.5033831434694\\
};

\end{axis}
\end{tikzpicture}%
\vspace*{-0.3in}
\caption{Passenger travel cost (including time and money) as a function of $K$.} 
\label{figure5}
\end{minipage}
\begin{minipage}[b]{0.005\linewidth}
\hfill 
\end{minipage}
\begin{minipage}[b]{0.32\linewidth}
\centering
%
%
\definecolor{mycolor1}{rgb}{0.00000,0.44700,0.74100}%
\definecolor{mycolor2}{rgb}{0.85000,0.32500,0.09800}%
\begin{tikzpicture}

\begin{axis}[%
width=1.7in,
height=0.8in,
at={(1.358in,0.0in)},
scale only axis,
xmin=0,
xmax=3000,
xtick={0,1000,2000,3000},
xticklabels={{0},{1e3},{2e3},{3e3}},
xlabel style={font=\color{white!15!black}},
xlabel={Parking Supply},
ymin=24,
ymax=30,
ylabel style={font=\color{white!15!black}},
ylabel={Driver Wage},
axis background/.style={fill=white},
legend style={at={(0.665,0)}, anchor=south west, legend cell align=left, align=left, draw=white!12!black}
]
\addplot [color=black, line width=1.0pt]
  table[row sep=crcr]{%
0	27.4787655378629\\
30.3030303030303	27.4787655378629\\
60.6060606060606	27.4787655378629\\
90.9090909090909	27.4787659481223\\
121.212121212121	27.4787655378629\\
151.515151515152	27.4787659481223\\
181.818181818182	27.4787659481223\\
212.121212121212	27.4787655378629\\
242.424242424242	27.4787655378629\\
272.727272727273	27.4787659481223\\
303.030303030303	27.4787655378629\\
333.333333333333	27.4787655378629\\
363.636363636364	27.4787659481223\\
393.939393939394	27.4787655378629\\
424.242424242424	27.4787655378629\\
454.545454545455	27.4787655378629\\
484.848484848485	27.4787655378629\\
515.151515151515	27.4787655378629\\
545.454545454545	27.4787655378629\\
575.757575757576	27.4787655378629\\
606.060606060606	27.4787659481223\\
636.363636363636	27.4787655378629\\
666.666666666667	27.4787655378629\\
696.969696969697	27.4787659481223\\
727.272727272727	27.4787655378629\\
757.575757575758	27.4787659481223\\
787.878787878788	27.4787655378629\\
818.181818181818	27.4787655378629\\
848.484848484848	27.4787659481223\\
878.787878787879	27.4787655378629\\
909.090909090909	27.4787655378629\\
939.393939393939	27.4787655378629\\
969.69696969697	27.4787659481223\\
1000	27.4787659481223\\
1030.30303030303	27.4787659481223\\
1060.60606060606	27.4787659481223\\
1090.90909090909	27.4787659481223\\
1121.21212121212	27.4787655378629\\
1151.51515151515	27.4787655378629\\
1181.81818181818	27.4787604779961\\
1212.12121212121	27.4786459472289\\
1242.42424242424	27.4810370935527\\
1272.72727272727	27.4872583592432\\
1303.0303030303	27.5245602931763\\
1333.33333333333	27.6176372501415\\
1363.63636363636	27.753695032322\\
1393.93939393939	27.9076616391917\\
1424.24242424242	28.0793077300118\\
1454.54545454545	28.250131594704\\
1484.84848484848	28.4262113348668\\
1515.15151515152	28.5833626837122\\
1545.45454545455	28.7578711586157\\
1575.75757575758	28.9016868693992\\
1606.06060606061	29.0032323433536\\
1636.36363636364	29.0807057824413\\
1666.66666666667	29.1223590422173\\
1696.9696969697	29.1223590422173\\
1727.27272727273	29.1223590422173\\
1757.57575757576	29.1223590422173\\
1787.87878787879	29.1223590422173\\
1818.18181818182	29.1223590422173\\
1848.48484848485	29.1223590422173\\
1878.78787878788	29.1223590422173\\
1909.09090909091	29.1223590422173\\
1939.39393939394	29.1223590422173\\
1969.69696969697	29.1223590422173\\
2000	29.1223590422173\\
2030.30303030303	29.1223590422173\\
2060.60606060606	29.1223590422173\\
2090.90909090909	29.1223590422173\\
2121.21212121212	29.1223590422173\\
2151.51515151515	29.1223590422173\\
2181.81818181818	29.1223590422173\\
2212.12121212121	29.1223590422173\\
2242.42424242424	29.1223590422173\\
2272.72727272727	29.1223590422173\\
2303.0303030303	29.1223590422173\\
2333.33333333333	29.1223590422173\\
2363.63636363636	29.1223590422173\\
2393.93939393939	29.1223590422173\\
2424.24242424242	29.1223590422173\\
2454.54545454545	29.1223590422173\\
2484.84848484848	29.1223590422173\\
2515.15151515152	29.1223590422173\\
2545.45454545455	29.1223590422173\\
2575.75757575758	29.1223590422173\\
2606.06060606061	29.1223590422173\\
2636.36363636364	29.1223590422173\\
2666.66666666667	29.1223590422173\\
2696.9696969697	29.1223590422173\\
2727.27272727273	29.1223590422173\\
2757.57575757576	29.1223590422173\\
2787.87878787879	29.1223590422173\\
2818.18181818182	29.1223590422173\\
2848.48484848485	29.1223590422173\\
2878.78787878788	29.1223590422173\\
2909.09090909091	29.1223590422173\\
2939.39393939394	29.1223590422173\\
2969.69696969697	29.1223590422173\\
3000	29.1223590422173\\
};
\addlegendentry{$w_n$}

\addplot [color=mycolor1, dashed, line width=1.0pt]
  table[row sep=crcr]{%
0	27.4787655378629\\
30.3030303030303	27.4051091270763\\
60.6060606060606	27.3331997424431\\
90.9090909090909	27.2622442762326\\
121.212121212121	27.1920171074901\\
151.515151515152	27.1224012848485\\
181.818181818182	27.0533203020587\\
212.121212121212	26.9847214948163\\
242.424242424242	26.9165664990289\\
272.727272727273	26.848824850824\\
303.030303030303	26.7814708407103\\
333.333333333333	26.7144856601121\\
363.636363636364	26.6478524964491\\
393.939393939394	26.5815557626164\\
424.242424242424	26.5155844240798\\
454.545454545455	26.4499273288384\\
484.848484848485	26.3845750578466\\
515.151515151515	26.3195192744972\\
545.454545454545	26.254752551414\\
575.757575757576	26.1902682333399\\
606.060606060606	26.1260607665426\\
636.363636363636	26.0621234119914\\
666.666666666667	25.9984525665728\\
696.969696969697	25.9350437509418\\
727.272727272727	25.8718915377983\\
757.575757575758	25.8089939515805\\
787.878787878788	25.7463457641863\\
818.181818181818	25.6839451447551\\
848.484848484848	25.6217891761565\\
878.787878787879	25.5598738113191\\
909.090909090909	25.4981979101124\\
939.393939393939	25.4367587178232\\
969.69696969697	25.3755545569265\\
1000	25.314582541423\\
1030.30303030303	25.2538412946593\\
1060.60606060606	25.193329119606\\
1090.90909090909	25.1330444434353\\
1121.21212121212	25.0729853465559\\
1151.51515151515	25.0131514032551\\
1181.81818181818	24.953535210697\\
1212.12121212121	24.8940179922594\\
1242.42424242424	24.837592433276\\
1272.72727272727	24.7859617337182\\
1303.0303030303	24.7710542553614\\
1333.33333333333	24.822409526991\\
1363.63636363636	24.9263085973361\\
1393.93939393939	25.0536843748763\\
1424.24242424242	25.2044054527916\\
1454.54545454545	25.3581683498277\\
1484.84848484848	25.5205542699662\\
1515.15151515152	25.6660937743715\\
1545.45454545455	25.8345118791199\\
1575.75757575758	25.9767632042527\\
1606.06060606061	26.0806421989769\\
1636.36363636364	26.1681119955397\\
1666.66666666667	26.2286065032254\\
1696.9696969697	26.2488879551416\\
1727.27272727273	26.2699089020158\\
1757.57575757576	26.2934240514495\\
1787.87878787879	26.3148937095221\\
1818.18181818182	26.3365814651395\\
1848.48484848485	26.3584860328969\\
1878.78787878788	26.3826915865037\\
1909.09090909091	26.404994562269\\
1939.39393939394	26.4275125018864\\
1969.69696969697	26.4522337659838\\
2000	26.4751504020848\\
2030.30303030303	26.4982808967921\\
2060.60606060606	26.52352582159\\
2090.90909090909	26.5470558290496\\
2121.21212121212	26.5707991876186\\
2151.51515151515	26.5965757200352\\
2181.81818181818	26.620719787647\\
2212.12121212121	26.6450767098327\\
2242.42424242424	26.6696463958804\\
2272.72727272727	26.6961524106475\\
2303.0303030303	26.7211248050371\\
2333.33333333333	26.7463103246603\\
2363.63636363636	26.7733659024205\\
2393.93939393939	26.798956631186\\
2424.24242424242	26.824760804129\\
2454.54545454545	26.8523742556645\\
2484.84848484848	26.8785865504191\\
2515.15151515152	26.9050132684209\\
2545.45454545455	26.9316543750451\\
2575.75757575758	26.9600311749893\\
2606.06060606061	26.9870846495578\\
2636.36363636364	27.0143539232763\\
2666.66666666667	27.04330815978\\
2696.9696969697	27.0709937142032\\
2727.27272727273	27.0988965374522\\
2757.57575757576	27.1284373312064\\
2787.87878787879	27.1567606383339\\
2818.18181818182	27.1853030052331\\
2848.48484848485	27.2140648201653\\
2878.78787878788	27.2444074479446\\
2909.09090909091	27.2735960688621\\
2939.39393939394	27.3030061297501\\
2969.69696969697	27.3339575395341\\
3000	27.363799791664\\
};
\addlegendentry{$w_g$}

\end{axis}
\end{tikzpicture}%
\vspace*{-0.3in}
\caption{Gross hourly wage ($w_g$) and net hourly wage of TNC drivers as a function of $K$.}
\label{figure6}
\end{minipage}
\end{figure*}

\begin{figure*}[bt]
\begin{minipage}[b]{0.32\linewidth}
\centering
%
%
\definecolor{mycolor1}{rgb}{0.00000,0.44700,0.74100}%
\begin{tikzpicture}

\begin{axis}[%
width=1.7in,
height=0.8in,
at={(1.358in,0.0in)},
scale only axis,
xmin=0,
xmax=3000,
xtick={0,1000,2000,3000},
xticklabels={{0},{1e3},{2e3},{3e3}},
xlabel style={font=\color{white!15!black}},
xlabel={Parking Supply},
ymin=4.4,
ymax=4.8,
ylabel style={font=\color{white!15!black}},
ylabel={Waiting Time (min)},
axis background/.style={fill=white},
legend style={legend cell align=left, align=left, draw=white!15!black}
]
\addplot [color=black, line width=1.0pt]
  table[row sep=crcr]{%
0	4.71298170449253\\
30.3030303030303	4.71298170449253\\
60.6060606060606	4.71298170449253\\
90.9090909090909	4.71298168399563\\
121.212121212121	4.71298170449253\\
151.515151515152	4.71298168399563\\
181.818181818182	4.71298168399563\\
212.121212121212	4.71298170449253\\
242.424242424242	4.71298170449253\\
272.727272727273	4.71298168399563\\
303.030303030303	4.71298170449253\\
333.333333333333	4.71298170449253\\
363.636363636364	4.71298168399563\\
393.939393939394	4.71298170449253\\
424.242424242424	4.71298170449253\\
454.545454545455	4.71298170449253\\
484.848484848485	4.71298170449253\\
515.151515151515	4.71298170449253\\
545.454545454545	4.71298170449253\\
575.757575757576	4.71298170449253\\
606.060606060606	4.71298168399563\\
636.363636363636	4.71298170449253\\
666.666666666667	4.71298170449253\\
696.969696969697	4.71298168399563\\
727.272727272727	4.71298170449253\\
757.575757575758	4.71298168399563\\
787.878787878788	4.71298170449253\\
818.181818181818	4.71298170449253\\
848.484848484848	4.71298168399563\\
878.787878787879	4.71298170449253\\
909.090909090909	4.71298170449253\\
939.393939393939	4.71298170449253\\
969.69696969697	4.71298168399563\\
1000	4.71298168399563\\
1030.30303030303	4.71298168399563\\
1060.60606060606	4.71298168399563\\
1090.90909090909	4.71298168399563\\
1121.21212121212	4.71298170449253\\
1151.51515151515	4.71298170449253\\
1181.81818181818	4.71298206523837\\
1212.12121212121	4.7129895453931\\
1242.42424242424	4.71266429687406\\
1272.72727272727	4.71137396360991\\
1303.0303030303	4.70467368417655\\
1333.33333333333	4.68874646909251\\
1363.63636363636	4.66518587434842\\
1393.93939393939	4.63819766470008\\
1424.24242424242	4.6092092638464\\
1454.54545454545	4.57916015032941\\
1484.84848484848	4.54961765368196\\
1515.15151515152	4.52202957802486\\
1545.45454545455	4.49454582303087\\
1575.75757575758	4.47058620371765\\
1606.06060606061	4.45302660365262\\
1636.36363636364	4.44099788069128\\
1666.66666666667	4.43525097124404\\
1696.9696969697	4.43494514176946\\
1727.27272727273	4.43490120200119\\
1757.57575757576	4.43494214977944\\
1787.87878787879	4.43492370160667\\
1818.18181818182	4.43490597826187\\
1848.48484848485	4.43488883212788\\
1878.78787878788	4.43492847797393\\
1909.09090909091	4.43491155953896\\
1939.39393939394	4.43489505729461\\
1969.69696969697	4.43493283844132\\
2000	4.43491657734811\\
2030.30303030303	4.4349007190098\\
2060.60606060606	4.43493672933584\\
2090.90909090909	4.43492111218107\\
2121.21212121212	4.43490597826187\\
2151.51515151515	4.43494036532326\\
2181.81818181818	4.43492539211609\\
2212.12121212121	4.43491075454517\\
2242.42424242424	4.43489651967771\\
2272.72727272727	4.43492928298187\\
2303.0303030303	4.43491512835253\\
2333.33333333333	4.43490144349697\\
2363.63636363636	4.43493283844132\\
2393.93939393939	4.43491924726209\\
2424.24242424242	4.43490597826187\\
2454.54545454545	4.43493615240916\\
2484.84848484848	4.4349229771035\\
2515.15151515152	4.43491019104989\\
2545.45454545455	4.4348976466533\\
2575.75757575758	4.43492653254071\\
2606.06060606061	4.43491406844041\\
2636.36363636364	4.43490193990515\\
2666.66666666667	4.43492976598695\\
2696.9696969697	4.43491771776276\\
2727.27272727273	4.43490597826187\\
2757.57575757576	4.43493283844132\\
2787.87878787879	4.43492111218107\\
2818.18181818182	4.43490970805419\\
2848.48484848485	4.434898532135\\
2878.78787878788	4.43492434560992\\
2909.09090909091	4.43491326344462\\
2939.39393939394	4.43490234239845\\
2969.69696969697	4.4349273375471\\
3000	4.43491657734811\\
};

\end{axis}
\end{tikzpicture}%
\vspace*{-0.3in}
\caption{Passenger waiting time (minute) as a function of $K$. }
\label{figure7}
\end{minipage}
\begin{minipage}[b]{0.005\linewidth}
\hfill
\end{minipage}
\begin{minipage}[b]{0.32\linewidth}
\centering
%
%
\definecolor{mycolor1}{rgb}{0.00000,0.44700,0.74100}%
\definecolor{mycolor2}{rgb}{0.85000,0.32500,0.09800}%
\definecolor{mycolor3}{rgb}{0.92900,0.69400,0.12500}%
\begin{tikzpicture}

\begin{axis}[%
width=1.7in,
height=0.8in,
at={(1.358in,0.0in)},
scale only axis,
xmin=0,
xmax=3000,
xtick={0,1000,2000,3000},
xticklabels={{0},{1e3},{2e3},{3e3}},
xlabel style={font=\color{white!15!black}},
xlabel={Parking Supply},
ymin=7,
ymax=16,
ylabel style={font=\color{white!15!black}},
ylabel={Price (\$/trip)},
axis background/.style={fill=white},
legend style={at={(0.67,0.4)}, anchor=south west, legend cell align=left, align=left, draw=white!12!black}
]
\addplot [color=black, line width=1.0pt]
  table[row sep=crcr]{%
0	14.6915045230725\\
30.3030303030303	14.6915045230725\\
60.6060606060606	14.6915045230725\\
90.9090909090909	14.6915043896545\\
121.212121212121	14.6915045230725\\
151.515151515152	14.6915043896545\\
181.818181818182	14.6915043896545\\
212.121212121212	14.6915045230725\\
242.424242424242	14.6915045230725\\
272.727272727273	14.6915043896545\\
303.030303030303	14.6915045230725\\
333.333333333333	14.6915045230725\\
363.636363636364	14.6915043896545\\
393.939393939394	14.6915045230725\\
424.242424242424	14.6915045230725\\
454.545454545455	14.6915045230725\\
484.848484848485	14.6915045230725\\
515.151515151515	14.6915045230725\\
545.454545454545	14.6915045230725\\
575.757575757576	14.6915045230725\\
606.060606060606	14.6915043896545\\
636.363636363636	14.6915045230725\\
666.666666666667	14.6915045230725\\
696.969696969697	14.6915043896545\\
727.272727272727	14.6915045230725\\
757.575757575758	14.6915043896545\\
787.878787878788	14.6915045230725\\
818.181818181818	14.6915045230725\\
848.484848484848	14.6915043896545\\
878.787878787879	14.6915045230725\\
909.090909090909	14.6915045230725\\
939.393939393939	14.6915045230725\\
969.69696969697	14.6915043896545\\
1000	14.6915043896545\\
1030.30303030303	14.6915043896545\\
1060.60606060606	14.6915043896545\\
1090.90909090909	14.6915043896545\\
1121.21212121212	14.6915045230725\\
1151.51515151515	14.6915045230725\\
1181.81818181818	14.6915056173152\\
1212.12121212121	14.6915338855918\\
1242.42424242424	14.6918071586203\\
1272.72727272727	14.6947870467694\\
1303.0303030303	14.7073870619517\\
1333.33333333333	14.7349587645383\\
1363.63636363636	14.7771546095221\\
1393.93939393939	14.8272791461968\\
1424.24242424242	14.8783253306922\\
1454.54545454545	14.9363144159512\\
1484.84848484848	14.9895667905675\\
1515.15151515152	15.0444166785155\\
1545.45454545455	15.0894987543425\\
1575.75757575758	15.1343647007747\\
1606.06060606061	15.1699293370377\\
1636.36363636364	15.1899091118904\\
1666.66666666667	15.1968207189578\\
1696.9696969697	15.1984785580655\\
1727.27272727273	15.1987167650232\\
1757.57575757576	15.1984947781413\\
1787.87878787879	15.1985947892315\\
1818.18181818182	15.1986908716596\\
1848.48484848485	15.1987838256098\\
1878.78787878788	15.1985688955482\\
1909.09090909091	15.1986606141601\\
1939.39393939394	15.1987500771719\\
1969.69696969697	15.1985452565927\\
2000	15.1986334114251\\
2030.30303030303	15.1987193834494\\
2060.60606060606	15.1985241633231\\
2090.90909090909	15.1986088270707\\
2121.21212121212	15.1986908716596\\
2151.51515151515	15.1985044519856\\
2181.81818181818	15.198585624621\\
2212.12121212121	15.1986649782282\\
2242.42424242424	15.1987421491816\\
2272.72727272727	15.1985645314376\\
2303.0303030303	15.198641266768\\
2333.33333333333	15.1987154558098\\
2363.63636363636	15.1985452565927\\
2393.93939393939	15.198618937212\\
2424.24242424242	15.1986908716596\\
2454.54545454545	15.1985272909486\\
2484.84848484848	15.1985987169191\\
2515.15151515152	15.1986680330748\\
2545.45454545455	15.1987360395334\\
2575.75757575758	15.1985794421409\\
2606.06060606061	15.1986470128019\\
2636.36363636364	15.1987127646484\\
2666.66666666667	15.1985619129703\\
2696.9696969697	15.1986272289749\\
2727.27272727273	15.1986908716596\\
2757.57575757576	15.1985452565927\\
2787.87878787879	15.1986088270707\\
2818.18181818182	15.1986706515139\\
2848.48484848485	15.1987312390929\\
2878.78787878788	15.1985912979523\\
2909.09090909091	15.1986513768761\\
2939.39393939394	15.1987105826251\\
2969.69696969697	15.1985750780349\\
3000	15.1986334114251\\
};
\addlegendentry{$p_f$}

\addplot [color=mycolor1, dashed, line width=1.0pt]
  table[row sep=crcr]{%
0	9.28305880006287\\
30.3030303030303	9.25817570291678\\
60.6060606060606	9.23388279773126\\
90.9090909090909	9.20991205429345\\
121.212121212121	9.18618754373558\\
151.515151515152	9.16266936807119\\
181.818181818182	9.13933197259958\\
212.121212121212	9.11615756517652\\
242.424242424242	9.09313299252821\\
272.727272727273	9.07024796387706\\
303.030303030303	9.04749408134506\\
333.333333333333	9.0248647034216\\
363.636363636364	9.00235415108662\\
393.939393939394	8.97995744392196\\
424.242424242424	8.95767056884709\\
454.545454545455	8.93548985352611\\
484.848484848485	8.91341211595459\\
515.151515151515	8.89143453980458\\
545.454545454545	8.86955461591061\\
575.757575757576	8.84777009594957\\
606.060606060606	8.82607901047813\\
636.363636363636	8.80447936260632\\
666.666666666667	8.78296965537234\\
696.969696969697	8.76154837621589\\
727.272727272727	8.74021396933938\\
757.575757575758	8.71896539754343\\
787.878787878788	8.69780126353811\\
818.181818181818	8.67672067247071\\
848.484848484848	8.65572263952507\\
878.787878787879	8.6348060718275\\
909.090909090909	8.6139703098767\\
939.393939393939	8.59321451450215\\
969.69696969697	8.5725380284251\\
1000	8.55194005802801\\
1030.30303030303	8.53142004745612\\
1060.60606060606	8.5109774234077\\
1090.90909090909	8.49061165453965\\
1121.21212121212	8.47032218184918\\
1151.51515151515	8.45010868233317\\
1181.81818181818	8.42996952691879\\
1212.12121212121	8.4098826879558\\
1242.42424242424	8.39090414129317\\
1272.72727272727	8.3749725699501\\
1303.0303030303	8.37600503452043\\
1333.33333333333	8.40644393187595\\
1363.63636363636	8.46220181247687\\
1393.93939393939	8.53067096679862\\
1424.24242424242	8.60794989195502\\
1454.54545454545	8.69139758325098\\
1484.84848484848	8.7757704789702\\
1515.15151515152	8.85673865553251\\
1545.45454545455	8.93993531852173\\
1575.75757575758	9.01536246283574\\
1606.06060606061	9.07281720737267\\
1636.36363636364	9.11505235883195\\
1666.66666666667	9.13997740138039\\
1696.9696969697	9.14815073542386\\
1727.27272727273	9.15563590677785\\
1757.57575757576	9.16368310510222\\
1787.87878787879	9.17123251454428\\
1818.18181818182	9.17885541631968\\
1848.48484848485	9.1865518934639\\
1878.78787878788	9.19484398015796\\
1909.09090909091	9.20267852118144\\
1939.39393939394	9.21058657198967\\
1969.69696969697	9.21906479217956\\
2000	9.22711095640802\\
2030.30303030303	9.23523028865648\\
2060.60606060606	9.24389709803689\\
2090.90909090909	9.25215482569094\\
2121.21212121212	9.26048524414758\\
2151.51515151515	9.26934289862773\\
2181.81818181818	9.27781242083642\\
2212.12121212121	9.28635499584914\\
2242.42424242424	9.29497034342948\\
2272.72727272727	9.30408777143938\\
2303.0303030303	9.31284321150216\\
2333.33333333333	9.32167129760752\\
2363.63636363636	9.33098494224525\\
2393.93939393939	9.33995394141007\\
2424.24242424242	9.34899623644661\\
2454.54545454545	9.35850847056041\\
2484.84848484848	9.36769269307859\\
2515.15151515152	9.37695029963775\\
2545.45454545455	9.38628182267239\\
2575.75757575758	9.3960644688076\\
2606.06060606061	9.40553944575419\\
2636.36363636364	9.41508847812544\\
2666.66666666667	9.42507597474704\\
2696.9696969697	9.43476982345801\\
2727.27272727273	9.44453833495135\\
2757.57575757576	9.45473352759301\\
2787.87878787879	9.46464856999888\\
2818.18181818182	9.47463884492261\\
2848.48484848485	9.48470483705189\\
2878.78787878788	9.49518301572868\\
2909.09090909091	9.50539743813538\\
2939.39393939394	9.51568852085992\\
2969.69696969697	9.5263816212434\\
3000	9.53682275764279\\
};
\addlegendentry{$p_d$}

\end{axis}

\end{tikzpicture}%
\vspace*{-0.3in}
\caption{Per-trip passenger fare and driver payment as a function of $K$.} 
\label{figure8}
\end{minipage}
\begin{minipage}[b]{0.005\linewidth}
\hfill 
\end{minipage}
\begin{minipage}[b]{0.32\linewidth}
\centering
%
%
\definecolor{mycolor1}{rgb}{0.00000,0.44700,0.74100}%
\begin{tikzpicture}

\begin{axis}[%
width=1.7in,
height=0.8in,
at={(1.358in,0.0in)},
scale only axis,
xmin=0,
xmax=3000,
xtick={0,1000,2000,3000},
xticklabels={{0},{1e3},{2e3},{3e3}},
xlabel style={font=\color{white!15!black}},
xlabel={Parking Supply},
ymin=48000,
ymax=60000,
ylabel style={font=\color{white!15!black}},
ylabel={Profit (\$/hour)},
axis background/.style={fill=white},
legend style={legend cell align=left, align=left, draw=white!15!black}
]
\addplot [color=black, line width=1.0pt]
  table[row sep=crcr]{%
0	48879.2560718392\\
30.3030303030303	49104.1390307435\\
60.6060606060606	49323.6880831155\\
90.9090909090909	49540.3259607069\\
121.212121212121	49754.7377141672\\
151.515151515152	49967.2854884361\\
181.818181818182	50178.1990529636\\
212.121212121212	50387.6391853109\\
242.424242424242	50595.7255820428\\
272.727272727273	50802.5512669559\\
303.030303030303	51008.1908522535\\
333.333333333333	51212.7056453966\\
363.636363636364	51416.1469907123\\
393.939393939394	51618.5585556383\\
424.242424242424	51819.9779521927\\
454.545454545455	52020.4379215732\\
484.848484848485	52219.9672213966\\
515.151515151515	52418.5913044497\\
545.454545454545	52616.3328475137\\
575.757575757576	52813.2121699866\\
606.060606060606	53009.2475699327\\
636.363636363636	53204.45559721\\
666.666666666667	53398.8512779223\\
696.969696969697	53592.4483007094\\
727.272727272727	53785.2591727526\\
757.575757575758	53977.2953514839\\
787.878787878788	54168.5673566066\\
818.181818181818	54359.0848660171\\
848.484848484848	54548.8567984511\\
878.787878787879	54737.8913851011\\
909.090909090909	54926.1962320024\\
939.393939393939	55113.7783746448\\
969.69696969697	55300.6443259955\\
1000	55486.8001189059\\
1030.30303030303	55672.2513437062\\
1060.60606060606	55857.0031816562\\
1090.90909090909	56041.0604348115\\
1121.21212121212	56224.4275527357\\
1151.51515151515	56407.1086536568\\
1181.81818181818	56589.1074301516\\
1212.12121212121	56770.4242220831\\
1242.42424242424	56951.0134357697\\
1272.72727272727	57130.3450583764\\
1303.0303030303	57305.7851267737\\
1333.33333333333	57471.052824327\\
1363.63636363636	57617.7125695609\\
1393.93939393939	57742.7587774559\\
1424.24242424242	57842.8607069436\\
1454.54545454545	57913.1440849298\\
1484.84848484848	57960.019982037\\
1515.15151515152	57981.1093726874\\
1545.45454545455	57977.6619558373\\
1575.75757575758	57945.0168265153\\
1606.06060606061	57899.4742176462\\
1636.36363636364	57840.1310928062\\
1666.66666666667	57765.9771243261\\
1696.9696969697	57696.8618919049\\
1727.27272727273	57626.7526804117\\
1757.57575757576	57548.8294200262\\
1787.87878787879	57477.3716979133\\
1818.18181818182	57405.19316956\\
1848.48484848485	57332.2972821483\\
1878.78787878788	57252.0722009038\\
1909.09090909091	57177.8535017736\\
1939.39393939394	57102.9226856602\\
1969.69696969697	57020.9731840683\\
2000	56944.7186708039\\
2030.30303030303	56867.7556396783\\
2060.60606060606	56784.0558231982\\
2090.90909090909	56705.7663803805\\
2121.21212121212	56626.7705582407\\
2151.51515151515	56541.2948433052\\
2181.81818181818	56460.9681812553\\
2212.12121212121	56379.9359643386\\
2242.42424242424	56298.1988705963\\
2272.72727272727	56210.2896459447\\
2303.0303030303	56127.214615401\\
2333.33333333333	56043.433883415\\
2363.63636363636	55953.6905140046\\
2393.93939393939	55868.5636078407\\
2424.24242424242	55782.7291142095\\
2454.54545454545	55691.1248068235\\
2484.84848484848	55603.9344578157\\
2515.15151515152	55516.033643664\\
2545.45454545455	55427.4216534412\\
2575.75757575758	55333.2731569792\\
2606.06060606061	55243.2911232624\\
2636.36363636364	55152.5937520549\\
2666.66666666667	55056.5202028964\\
2696.9696969697	54964.4397833728\\
2727.27272727273	54871.6389920327\\
2757.57575757576	54773.6106042739\\
2787.87878787879	54679.4124189739\\
2818.18181818182	54584.4879851336\\
2848.48484848485	54488.8354868108\\
2878.78787878788	54388.1363528001\\
2909.09090909091	54291.0658991879\\
2939.39393939394	54193.2603854272\\
2969.69696969697	54090.5333758255\\
3000	53991.2923840913\\
};

\end{axis}
\end{tikzpicture}%
\vspace*{-0.3in}
\caption{Hourly platform profit (\$/hour) as a function of $K$.}
\label{figure9}
\end{minipage}
\end{figure*}

Drivers choose their work based on the net hourly wage. Under a logit model, the driver supply function can be written as:
\begin{equation}
N=N_0\dfrac{e^{\eta w_n}}{e^{\eta w_n}+e^{\eta w_0}},
\end{equation}
where $N_0$ is the number of potential drivers, and $\eta$ and $w_0$ are the model parameters. 
Garage operators decide whether to offer the parking space to the TNC platform based on the earnings offered by the platform. Assume that each parking space has a reservation earning, and the garage operator only offers the parking space if the platform earning $e$ exceeds the reservation earning. Different parking spaces have distinct reservation earnings. We assume that the distribution of the reservation earning is captured by a log-normal distribution. In this case, the parking supply function $F_g$ can be written as:
\begin{equation}
   F_g(p_gr)=K_0\left(\dfrac{1}{2}+\dfrac{1}{2}\text{erf}\left(\dfrac{\text{ln}(p_gr)-u_0}{\sqrt{2}\sigma}\right)\right)
\end{equation}
where $K_0$ is all potential {\em idle} parking space calculated based on the arrival and departure pattern of long-term regular parking demand, $\text{erf}(\cdot)$ denotes the error function, and $u_0$ and $\sigma$ are the model parameters. 

To setup the numerical study, we need to specify the values of the following parameters:
\begin{equation}
    \Theta=\{\lambda_0, N_0, K_0, M, \alpha, \epsilon, c_0, \eta, w_0, \sigma, u_0, l\}.
\end{equation}
To this end, we consider the following profit maximization problem without parking services:
\begin{align}
\label{TNCprofit_max_noparking}
&\max_{p_f, w_g} \lambda p_f-w_gN \\
& 
\begin{cases}
\label{marekt_constraints_noparking}
\lambda=F_p\bigg( \alpha t_w(N-\lambda/\mu)+p_f\bigg), \\
N=F_d\left( w_g \right)  
\end{cases}
\end{align}
and we set the parameters values $\Theta$ so that the optimal solution to (\ref{TNCprofit_max_noparking}) matches the San Francisco TNC market data. The values of these parameters are summarized below:
\begin{align*}
    \lambda_0=944/\text{min}, N_0=&10,000, K_0=10,000,  M = 174.7, \\
    \alpha=3, \epsilon=0.155, c_0&=15.48, \eta=0.144, w_0=32.41, \\
    \sigma=0.6, &u_0=1.1, l=\$8/\text{hour}
\end{align*}
The detailed justification of these parameters values can be found in Appendix 1 of \cite{li2020congestion}.

To understand impacts of parking on the TNC business model, we  solve (\ref{TNCprofit_max}) in two steps: first we fix $K$ and derive other endogenous variables as a function of $K$, which reflect how parking provision affects the surplus of passenger, drivers, garage operators, and the TNC platform; then we vary $K$ to derive the optimal parking supply that maximizes the platform profit. Note that $K$ is {\em endogenous} and the platform eventually chooses $K$ that maximizes its profit. 

As the number of parking spaces $K$ increases, the average occupancy of each garage decreases. Therefore, it is easier for each idle TNC vehicle to find a parking space. This indicates that the occupancy of parking garages is a decreasing function of $K$, and the ratio of parked TNC vehicles out of all {\em idle} TNC vehicles is an increasing function of $K$. The average occupancy of the TNC parking spaces is shown in Figure 2. The ratio of parked TNC vehicles is shown in Figure 3. The number of passengers and drivers as are shown in Figure 4 and Figure 5, respectively. The passenger travel cost is shown in Figure 6. The driver gross wage and net wage are shown in Figure 7, and passenger waiting time is shown in Figure 8. The passenger ride fare, per-trip driver payment, and per hour parking rate is shown in Figure 9. The platform profit is shown in Figure 10. All variables are represented as functions of $K$.

Based on Figure 2-10, the optimal solution as a function of $K$ clearly has three distinct regimes:
\begin{itemize}
    \item When $K\leq 1272$, all parking spaces are fully occupied. Passengers and drivers are unaffected. Driver savings from parking the TNC vehicle are reaped by the TNC platform, who pays a reduced gross wage and enjoys an increased platform profit. The number of TNC vehicles on road is reduced, and therefore traffic congestion is reduced. 
    \item When $1272\leq K \leq 1607$, both passengers and drivers benefit from the parking services. The passenger total travel cost reduces, the driver net wage increases, more passenger are delivered, and more drivers are hired. The number of TNC vehicle on the road network slightly increases, but is significantly smaller than at $K=0$.  The platform profit is maximized at $K=1515$.
    \item When $K\geq 1607$, all idle TNC vehicles are parked. Passenger and drivers are unaffected. The cost due to over-provision of parking spaces are covered the TNC platform, who pays an increased gross wage and earns a reduced profit. Traffic congestion is unaffected. 
\end{itemize}
By comparing the proposed business model with  (\ref{TNCprofit_max_noparking}), our numerical study clearly indicates that at the optimal solution, the proposed business model benefits both passengers, drivers, garage operators, the TNC platform. In the meanwhile, it reduced the number of  TNC vehicles on the road network, which helps to alleviate the traffic congestion. 

\subsection{Analysis}
We emphasize that many of the aforementioned insights hold for a large range of model parameters. For notation convenience, let $\lambda^*(K)$,  $N^*(K)$, $c^*(K)$, $w_n^*(K)$ and  $\text{Profit}(K)$ denote the optimal to (\ref{TNCprofit_max}) as a function of $K$, where $\lambda^*(K)$ denotes the passenger arrival rate,  $N^*(K)$ denotes the number of drivers, $c^*(K)$ denotes the passenger travel cost, $w_n^*(K)$ denotes the driver net wage, and  $\text{Profit}(K)$ represents the platform profit.  Let $K^*$ denote the profit-maximizing parking provision. We introduce the following proposition:
\vspace{0.1cm}
\begin{proposition}
Assume that the profit maximization problem (\ref{TNCprofit_max}) with exogenous $K$ has a unique solution, then there exists $K_1<K_2<K_0$ such that\\
(i) $\lambda^*(K)=\lambda^*(0)$, $N^*(K)=N^*(0)$,  $\forall 0\leq K\leq K_1$, \\ (ii) $\lambda^*(K)=\lambda^*(K_2)$, $N^*(K)=N^*(K_2)$ , $\forall K_2\leq K\leq \bar{K}$, \\
(iii) $\lambda^*(K^*)>\lambda^*(0)$ and $c^*(K^*)<c^*(0)$ \\
(iv) $N^*(K^*)>N^*(0)$ and  $w_n^*(K^*)>w_n^*(0)$, \\
(v) $\text{Profit}^*(K^*)>\text{Profit}(0)$.
\end{proposition}
\vspace{0.1cm}
The proof of Proposition 1 is deferred to the journal version due to space limit. The assumption of solution uniqueness is numerically validated and holds for all parameters within regime of practical interest. In Proposition 1,  (i) and (ii) represents the insight for regime 1 and regime 3, and (iii)-(v) shows that at the optimal solution $K^*$,  parking services benefit passengers, drivers and the TNC platform. 

In the case study, the optimal number of parking space is $K^*=1515$. This only accounts for less than 1\% of the 166,500 parking spots in the parking garages and commercial lots of San Francisco (note that $K_0$ denotes the idle parking spots  out of all spots).  At this optimal solution, compared to the case of $K=0$,  passenger arrival rate increases by $3.7\%$ (from $150.6$/min to $156.2$/min), the  number of drivers increases by 5.9\% (from 3053 to 3234), the number of on-road drivers reduces by 43\% (from 3053 to 1748), the passenger travel cost reduces by 0.8\% (from 28.83/trip to 28.61/trip), the driver net hour wage increases by 4\% (from \$27.48/hour to \$28.58/hour), and the platform profit increases by 18.6\% (from \$48,879/hour to \$57,981/hour).

\section{conclusion}
This paper investigated a novel business model that integrates parking services in the ride-hailing platform to reduce cruising congestion of TNC vehicles. The business model consists of a TNC platform, a group of passengers and drivers, and a number of parking garages. The platform provides app-based user interfaces that connect the passengers, drivers, and the parking garages. It matches passengers to drivers for ride-hailing services, and matches vacant TNC vehicles to unused parking garages to reduce cruising costs. The central thesis is that by enabling the sharing of idle parking spaces between garage operators and the TNC platform, the proposed business model will unlock the potential of the unused parking resources to benefit passengers, drivers, garage operators, and the TNC platform, and in the meanwhile curb cruising congestion without affecting regular parking demand. 

\bibliographystyle{unsrt}
\bibliography{TransactiveControl}

\begin{thebibliography}{10}

\bibitem{Uber_S1}
Uber.
\newblock {Form S1 Registration Statement}.
\newblock U. S. Securities and Exchange Commission, 2019.

\bibitem{cohen2016using}
P.~Cohen, R.~Hahn, J.~Hall, S.~Levitt, and R.~Metcalfe.
\newblock Using big data to estimate consumer surplus: The case of uber.
\newblock Technical report, National Bureau of Economic Research, 2016.

\bibitem{schaller2017empty}
B.~Schaller.
\newblock {Empty Seats, Full Streets. Fixing Manhattan Traffic Problem}.
\newblock Schaller Consulting, December 2017.

\bibitem{NYC2019improving}
NYCTLC.
\newblock {Improving efficiency and managing growth in New York's for-hire
  vehicle sector}.
\newblock New York City Taxi and Limousine Commission, 2019.

\bibitem{li2019regulating}
S.~Li, H.~Tavafoghi, K.~Poolla, and P.~Varaiya.
\newblock Regulating tncs: Should uber and lyft set their own rules?
\newblock {\em Transportation Research Part B: Methodological}, 2019.

\bibitem{li2020congestion}
S.~Li, K.~Poolla, and P.~Varaiya.
\newblock Impacts of congestion charge and minimum wage on tncs: a case study
  for san francsico.
\newblock {\em arXiv preprint arXiv:2003.02550}, 2020.

\bibitem{NYC2019crusing}
CNBC.
\newblock In win for uber and lyft, judge strikes down new york city’s
  cruising cap.
\newblock Reuters, 2019.
\newblock
  {https://www.cnbc.com/2019/12/23/in-win-for-uber-judge-strikes-down-new-york-citys-cruising-cap.html}.

\bibitem{xu2017optimal}
Z.~Xu, Y.~Yin, and L.~Zha.
\newblock Optimal parking provision for ride-sourcing services.
\newblock {\em Transportation Research Part B: Methodological}, 105:559--578,
  2017.

\bibitem{ruch2019many}
C.~Ruch, S.~H{\"o}rl, R.~Ehrler, M.~Balac, and E.~Frazzoli.
\newblock How many parking spaces does a mobility-on-demand system require?
\newblock 2019.

\bibitem{kondor2018large}
D.~Kondor, P.~Santi, K.~Basak, X.~Zhang, and C.~Ratti.
\newblock Large-scale estimation of parking requirements for autonomous
  mobility on demand systems.
\newblock {\em arXiv preprint arXiv:1808.05935}, 2018.

\bibitem{lam2017coordinated}
A.~Y. Lam, J.~James, Y.~Hou, and V.~Li.
\newblock Coordinated autonomous vehicle parking for vehicle-to-grid services:
  Formulation and distributed algorithm.
\newblock {\em IEEE Transactions on Smart Grid}, 9(5):4356--4366, 2017.

\bibitem{jian2019unlock}
S.~Jian, W.~Liu, X.~Wang, H.~Yang, and S.~T. Waller.
\newblock Unlock the sharing economy: Integrating carsharing and parking
  sharing services.
\newblock Technical report, 2019.

\bibitem{SFpark}
SFpark.
\newblock {http://sfpark.org/}.

\bibitem{korolko2018dynamic}
N.~Korolko, D.~Woodard, C.~Yan, and H.~Zhu.
\newblock Dynamic pricing and matching in ride-hailing platforms.
\newblock {\em Available at SSRN}, 2018.

\end{thebibliography}
\end{document}